\documentclass[a4paper,11pt]{article}
\usepackage[active]{srcltx}
\usepackage{amsmath}\usepackage{amssymb}\usepackage{amscd}\usepackage{amsthm}
\usepackage{empheq}

\newtheorem{defi}{Definition}\newtheorem{theo}[defi]{Theorem}\newtheorem{prop}[defi]{Proposition}\newtheorem{lemm}[defi]{Lemma}      
  
\newtheorem*{theo-intro}{Theorem}

\newcounter{pic}

\newcommand{\C}{\mathbb{C}}
\newcommand{\tG}{\widetilde{G}}
\newcommand{\tj}{\widetilde{j}}

\newcommand{\lam}{\underline{\lambda}}
\newcommand{\thm}{\underline{\theta}}
\newcommand{\mum}{\underline{\mu}}
\newcommand{\kam}{\underline{\kappa}}
\newcommand{\vm}{\underline{v}}
\newcommand{\cc}{\mathrm{c}}
\newcommand{\pos}{\mathrm{p}}
\newcommand{\di}{\mathrm{d}}

\begin{document}

$\ $

\begin{center}
\bigskip\medskip

{\Large\bf Fusion procedure for wreath products of finite groups by the symmetric group}

\vspace{1cm}
{\large L. Poulain d'Andecy}

\vskip 0.8cm

{\small Mathematics Laboratory of Versailles (CNRS UMR 8100), Versailles-Saint-Quentin University,\\
45 Avenue des Etats-Unis 78000 Versailles, France}

\end{center}

\vskip 0.8cm
\begin{abstract}
Let $G$ be a finite group. A complete system of pairwise orthogonal idempotents is constructed for the wreath product of $G$ by the symmetric group by means of a fusion procedure, that is by consecutive evaluations of a rational function with values in the group ring. This complete system of idempotents is indexed by standard Young multi-tableaux. Associated to the wreath product of $G$ by the symmetric group, a Baxterized form for the Artin generators of the symmetric group is defined and appears in the rational function used in the fusion procedure.
\end{abstract}

\section{{\hspace{-0.55cm}.\hspace{0.55cm}}Introduction}

The fusion procedure for the symmetric group $S_n$ originates in \cite{Ju} and has then been developed in the more general situation of Hecke algebras \cite{Ch}. The procedure allows to express a complete set of primitive idempotents, indexed by standard Young tableaux, of the group ring $\C S_n$ via a certain limiting process on a rational function in several variables with values in $\C S_n$.

\vskip .2cm
In \cite{Mo}, an alternative approach of the fusion procedure for the symmetric group has been proposed, and is based on the existence of a maximal commutative set in $\C S_n$ formed by the Jucys--Murphy elements. More precisely, this approach relies on an expression for the primitive idempotents in terms of the Jucys--Murphy elements \cite{Ju2,Mu} and transforms this expression into a fusion formula. The evaluations of the variables of the rational function are now consecutive. This version of the fusion procedure has then been generalized to the Hecke algebras of type A \cite{IMOs}, to the Brauer algebras \cite{IM,IMOg}, to the Birman--Wenzl--Murakami algebras \cite{IMOg2}, to the complex reflection groups of type $G(m,1,n)$ \cite{OPdA} and to the Ariki--Koike algebras \cite{OPdA2}.

\vskip .2cm
Let $G$ be a finite group. The aim of this article is to generalize the fusion procedure, in the spirit of \cite{Mo}, for the wreath product $\tG_n$ of the group $G$ by the symmetric group $S_n$. The representation theory of the group $\tG_n$ is well-known, see \emph{e.g.} \cite{JK} or \cite{Mac}, and the irreducible representations are labelled by multi-partitions of size $n$. Analogues of the Jucys--Murphy elements for the group $\tG_n$ have been introduced in \cite{Pus,Wan} and used for an inductive approach to the representation theory of the chain of groups $\tG_n$ \cite{Pus}.

\vskip .2cm
The irreducible representation $V_{\lam}$ of $\tG_n$ corresponding to the multi-partition $\lam$ admits the following decomposition, as a vector space,
 \[V_{\lam}=\bigoplus W_{{\cal{T}}}\ ,\]
 where the direct sum is over the set of standard Young multi-tableaux of shape $\lam$. The subspaces $W_{{\cal{T}}}$ are common eigenspaces for the Jucys--Murphy elements $\tj_1,\dots,\tj_n$ of the group $\tG_n$. Let $E_{{\cal{T}}}$ be the idempotent of $\C\tG_n$ associated to the subspace $W_{{\cal{T}}}$. The idempotents $E_{{\cal{T}}}$, where ${\cal{T}}$ runs through the set of standard multi-tableaux of size $n$, form together a complete system of pairwise orthogonal idempotents of $\C\tG_n$.

\vskip .2cm
We state here the main result of the article.
\begin{theo-intro}
The idempotent $E_{{\cal{T}}}$ of $\C\tG_n$ corresponding to the standard multi-tableau ${\cal{T}}$ can be obtained by the following consecutive evaluations
\begin{equation}\label{eq-idem-fin1-intro}
E_{{\cal{T}}}=\frac{1}{\mathsf{F}^G_{\lam}\,\mathsf{F}^{\phantom{G}}_{\lam}}\,\Phi(u_1,\dots,u_n,\underline{v_1},\dots,\underline{v_n})\Bigr\rvert_{v^{(\alpha)}_i=\xi_{\pos_i}^{(\alpha)},\!\textrm{\scriptsize{$\begin{array}{l}i=1,\dots,n\\[0.1em]\alpha=1,\dots,m\end{array}$}}}\ \Bigr\rvert_{u_1=\cc^G_1}\dots\ \dots\ \Bigr\rvert_{u_n=\cc^G_n}\ .
\end{equation}
\end{theo-intro}
\noindent The function $\Phi$ appearing in (\ref{eq-idem-fin1-intro}) is a rational function in several variables with values in the group ring $\C\tG_n$ and $\mathsf{F}^G_{\lam}$, $\mathsf{F}^{\phantom{G}}_{\lam}$ are complex numbers.

\vskip .2cm
The eigenvalues of the Jucys--Murphy elements are not sufficient to distinguish between the different subspaces $W_{{\cal{T}}}$. Therefore the commutative family used here for the fusion procedure is formed by the Jucys--Murphy elements $\tj_1,\dots,\tj_n$, together with an additional set of elements $\{g_i^{(\alpha)}\}$ of $\C \tG_n$. The elements $g_i^{(\alpha)}$ are images in $\C \tG_n$ of a set of elements $\{g^{(\alpha)}\}$ which linearly span the center of $\C G$ - the indices $i=1,\dots,n$ indicate in which copies of $\C G$ in $\C \tG_n$ they belong  (see Sections \ref{sec-def} and \ref{sec-wr} for precise definitions).

The set of variables is split into two parts: the variables $v^{(\alpha)}_i$ are first evaluated simultaneously at complex numbers $\xi_{\nu}^{(\alpha)}$ which are eigenvalues of elements $g_i^{(\alpha)}$. These variables correspond to the positions of the nodes of a multi-partition (their places in the $m$-tuple). Then the variables $u_1,\dots,u_n$ are consecutively evaluated at the eigenvalues of the Jucys--Murphy elements $\tj_1,\dots,\tj_n$ and are related to the classical contents of the nodes of a multi-partition. 

\vskip .2cm
The rational function $\Phi$ can be written as the product of two functions, the first one containing the variables $u_1,\dots,u_n$ and the second one containing the variables $v^{(\alpha)}_i$. As in the fusion procedure for the symmetric groups \cite{Mo}, the function related to the contents (that is, containing the variables $u_1,\dots,u_n$) is build up from a ``Baxterized'' form of the generators of $S_n$ inside $\C \tG_n$. However, the Baxterized form used here is a non-trivial generalization, associated to the finite group $G$, of the usual Baxterization for the symmetric groups.

\vskip .2cm
The coefficient $(\mathsf{F}^G_{\lam}\,\mathsf{F}^{\phantom{G}}_{\lam})^{-1}$ appearing in (\ref{eq-idem-fin1-intro}) only depends on the multi-partition $\lam$ and not on the standard multi-tableau ${\cal{T}}$ of shape $\lam$. The element $\mathsf{F}^{\phantom{G}}_{\lam}$ is the product of the hook lengths of the nodes in the multi-partition $\lam$ and is independent of the group $G$. The additional factor $\mathsf{F}^G_{\lam}$ depends on the numbers $\xi_{\nu}^{(\alpha)}$, and in turn on the choice of the set $\{g^{(\alpha)}\}$ of central elements of $\C G$.

\vskip .2cm
In general the subspaces $W_{\cal{T}}$ are not one-dimensional, and therefore the idempotent $E_{{\cal{T}}}$ of $\C\tG_n$ is not in general primitive. In fact, the subspaces $W_{\cal{T}}$ are all one-dimensional if and only if $G$ is an Abelian finite group. In this situation the idempotents $E_{{\cal{T}}}$, where ${\cal{T}}$ runs through the set of standard multi-tableaux of size $n$, form together a complete system of primitive pairwise orthogonal idempotents of $\C\tG_n$.

\vskip .2cm
For an Abelian finite group $G$, we provide a simplified fusion procedure for the wreath product $\tG_n$. This simplified version is obtained by replacing, in the function $\Phi$ appearing in (\ref{eq-idem-fin1-intro}), the set $\{g^{(\alpha)}\}$ of central elements of $\C G$ by a smaller subset.

If $G$ is the cyclic group of order $m$ then the group $\tG_n$ is isomorphic to the complex reflection group of type $G(m,1,n)$. A fusion procedure for the complex reflection group of type $G(m,1,n)$ has been given in \cite{OPdA}. The fusion procedure presented here, in its simplified version, is slightly different than the one in \cite{OPdA}, namely the Baxterized forms are different. The Baxterized form used here, when $G$ is the cyclic group of order 2, has been used in \cite{Ch2} for a fusion procedure for the Coxeter groups of type B (that is, the complex reflection groups of type $G(2,1,n)$).

\vskip .4cm
The paper is organized as follows.
Sections \ref{sec-def} and \ref{sec-wr} contain definitions and notations about the finite group $G$ and the wreath product $\tG_n$. 
The Jucys--Murphy elements for the group $\tG_n$ are defined in Section \ref{sec-JM}. 
In Section \ref{sec-bax} we introduce the Baxterized form, associated to the wreath product $\tG_n$, for the Artin generators of the symmetric group and prove that these Baxterized elements satisfy the Yang--Baxter equation with spectral parameters. 
In Section \ref{sec-tab}, we recall standard results on the representation theory of the groups $\tG_n$ and give the formula for the idempotents $E_{{\cal{T}}}$ in terms of the elements $\tj_1,\dots,\tj_n$ and $g_i^{(\alpha)}$.
The main result of the article, which gives the fusion procedure for the group $\tG_n$, is proved in Section \ref{sec-fus}. In Section \ref{sec-abe}, we consider the case of an Abelian finite group $G$ and provide a simplified fusion formula in this situation.

\paragraph{Notation.}$\ $

\vskip .2cm
\noindent We denote by $\C H$ the group ring over the complex numbers of a finite group $H$. For a vector space $V$, we denote by $\text{Id}_V$ the identity operator on $V$. Symbols $\vm$ and $\underline{v_i}$, for an integer $i$, stand for $m$-tuples of variables, namely $\vm:=(v^{(1)},\dots,v^{(m)})$ and $\underline{v_i}:=(v^{(1)}_i,\dots,v^{(m)}_i)$;

\section{{\hspace{-0.55cm}.\hspace{0.55cm}}Definitions}\label{sec-def}

Let $G$ be a finite group and let $\{C_1,\dots,C_m\}$ be the set of all its conjugacy classes. We denote by $g^{(1)},\dots,g^{(m)}$ the following central elements of the group ring $\C G$:
\begin{equation}\label{def-g}g^{(\alpha)}:=\frac{1}{|C_{\alpha}|}\sum\limits_{g\in C_{\alpha}}g\ ,\ \ \ \alpha=1,\dots,m.
\end{equation}

Let $\rho_1,\dots,\rho_m$ be the pairwise non-isomorphic irreducible representations of $G$ and $W_1,\dots,W_m$ be the corresponding representation spaces of dimensions, respectively, $\di_1$, $\dots$, $\di_m$. Denote by $\chi_1,\dots,\chi_m$ the associated irreducible characters. We define complex numbers $\xi_{\nu}^{(\alpha)}$, $\alpha,\nu=1,\dots,m$, by:
\begin{equation}\label{def-p}\xi_{\nu}^{(\alpha)}:=\frac{1}{\di_{\nu}}\chi_{\nu}(g^{(\alpha)})\ ,\ \ \ \alpha,\nu=1,\dots,m.
\end{equation}
The central elements $g^{(\alpha)}$ act in the irreducible representations of $G$ as multiples of the identity operators, namely:
$$\rho_{\nu}(g^{(\alpha)})=\xi_{\nu}^{(\alpha)}\cdot\mathrm{Id}_{W_{\nu}}\ \ \ \text{for $\alpha,\nu=1,\dots,m$\,.}$$
It is a standard fact that the elements $g^{(\alpha)}$, $\alpha=1,\dots,m$, span the center of $\C G$ and therefore, if $\nu\neq \nu'$, there exists some $\alpha\in\{1,\dots,m\}$ such that $\xi_{\nu\phantom{'}}^{(\alpha)}\neq \xi_{\nu'}^{(\alpha)}$; that is, the eigenvalues of the elements $g^{(\alpha)}$, $\alpha=1,\dots,m$, distinguish between the irreducible representations of $G$.

\paragraph{Functions $g^{(\alpha)}(v)$.} For $\alpha=1,\dots,m$, we define $S^{(\alpha)}$ to be the set formed by the pairwise different numbers among the $\xi_i^{(\alpha)}$, $i=1,\dots,m$, that is:
$$S^{(\alpha)}:=\{\xi_{i_1}^{(\alpha)},\xi_{i_2}^{(\alpha)},\dots,\xi_{i_{k_{\alpha}}}^{(\alpha)}\}\ \ \ \text{with $1=i_1<i_2<\dots <i_{k_{\alpha}}\leq m$ ,}$$
such that the numbers $\xi_{i_1}^{(\alpha)},\xi_{i_2}^{(\alpha)},\dots,\xi_{i_{k_{\alpha}}}^{(\alpha)}$ are pairwise different and, for any $i\in\{1,\dots,m\}$, we have $\xi_i^{(\alpha)}=\xi_{i_a}^{(\alpha)}$ for some $a\in\{1,\dots,k_{\alpha}\}$. By construction, we have in $\C G$ that:
\begin{equation}\label{char-g}\prod_{\xi^{(\alpha)}\in S^{(\alpha)}}\left(g^{(\alpha)}-\xi^{(\alpha)}\right)=0\ \ \ \ \text{for $\alpha=1,\dots,m$.}\end{equation}
 
We define the following rational functions in $v$ with values in $\C G$:
\begin{equation}\label{bax-g}
g^{(\alpha)}(v):=\frac{\prod_{\xi^{(\alpha)}\in S^{(\alpha)}}(v-\xi^{(\alpha)})}{v-g^{(\alpha)}}\ \ \ \ \text{for $\alpha=1,\dots,m$.}
\end{equation}
The rational functions $g^{(\alpha)}(v)$ can be rewritten as polynomial functions in $v$ as follows. Fix $\alpha\in\{1,\dots,m\}$ and let $a^{(\alpha)}_0,a^{(\alpha)}_1,\dots,a^{(\alpha)}_{k_{\alpha}}$ be the complex numbers defined by 
\[\prod_{\xi^{(\alpha)}\in S^{(\alpha)}}(X-\xi^{(\alpha)})=a^{(\alpha)}_0+a^{(\alpha)}_1X+\dots+a^{(\alpha)}_{k_{\alpha}}X^{k_{\alpha}}\ ,\]
where $X$ is an indeterminate. Define the polynomials $\mathfrak{a}^{(\alpha)}_i(v)$, $i=0,\dots,k_{\alpha}$, in $v$ by 
\[\mathfrak{a}^{(\alpha)}_i(v)=a^{(\alpha)}_i+a^{(\alpha)}_{i+1}\,v+\dots+a^{(\alpha)}_{k_{\alpha}}v^{k_{\alpha}-i}\ \quad\textrm{for $i=0,\dots,k_{\alpha}$.}\]
Then we have that:
\begin{equation}\label{bax-g2}
g^{(\alpha)}(v)=\sum_{i=0}^{k_{\alpha}-1}\mathfrak{a}^{(\alpha)}_{i+1}\!(v)\bigl(g^{(\alpha)}\bigr)^i=\mathfrak{a}^{(\alpha)}_1(v)+\mathfrak{a}^{(\alpha)}_2(v)g^{(\alpha)}+\dots+\mathfrak{a}^{(\alpha)}_{k_{\alpha}}(v)\bigl(g^{(\alpha)}\bigr)^{k_{\alpha}-1}\ .
\end{equation}
Indeed one can directly verify that $\displaystyle (v-g^{(\alpha)})\sum\limits_{i=0}^{k_{\alpha}-1}\mathfrak{a}^{(\alpha)}_{i+1}\!(v)\bigl(g^{(\alpha)}\bigr)^i=\prod_{\xi^{(\alpha)}\in S^{(\alpha)}}(v-\xi^{(\alpha)})\ $ (in the verification, it is useful to use the recursive relation $\mathfrak{a}^{(\alpha)}_{i+1}(v)=v^{-1}(\mathfrak{a}^{(\alpha)}_i(v)-a^{(\alpha)}_i)$ for $i=0,\dots,k_{\alpha}-1$, together with the initial condition $\mathfrak{a}^{(\alpha)}_0(v)=a^{(\alpha)}_0+a^{(\alpha)}_1\,v+\dots+a_{k_{\alpha}}^{(\alpha)}v^{k_{\alpha}}$).

\vskip .2cm
\noindent \textbf{Examples.} 

$-$ If $S^{(\alpha)}=\{\xi^{(\alpha)}\}$ then we have $g^{(\alpha)}(v)=1$;

$-$ if $S^{(\alpha)}=\{\xi^{(\alpha)}_1,\xi^{(\alpha)}_2\}$ then we have $g^{(\alpha)}(v)=g^{(\alpha)}+v-\xi^{(\alpha)}_1-\xi^{(\alpha)}_2$;  

$-$ if $S^{(\alpha)}=\{\xi^{(\alpha)}_1,\xi^{(\alpha)}_2,\xi^{(\alpha)}_3\}$ then we have\\
$\ $\hspace{0.3cm} $g^{(\alpha)}(v)=\bigl(g^{(\alpha)}\bigr)^2+(v-\xi^{(\alpha)}_1-\xi^{(\alpha)}_2-\xi^{(\alpha)}_3)g^{(\alpha)}+v^2-v(\xi^{(\alpha)}_1+\xi^{(\alpha)}_2+\xi^{(\alpha)}_3)+ \xi^{(\alpha)}_1\xi^{(\alpha)}_2+\xi^{(\alpha)}_1\xi^{(\alpha)}_3+\xi^{(\alpha)}_2\xi^{(\alpha)}_3\ $.

\section{{\hspace{-0.55cm}.\hspace{0.55cm}}Wreath product}\label{sec-wr}

\paragraph{Definitions.} Let $G^n:=G\times\dots\times G$, the Cartesian product of $n$ copies of $G$. The symmetric group $S_n$ on $n$ letters acts on $G^n$ by permuting the $n$ copies of $G$. We denote the action by $\sigma(a)$, $\sigma\in S_n$ and $a\in G^n$. The wreath product $\tG_n$ of the group $G$ by the symmetric group $S_n$ is the semi-direct product $G^n\rtimes S_n$ defined by this action. The group $\tG_n$ consists of the elements $(a,\sigma)$, $a\in G^n$ and $\sigma\in S_n$, with multiplication given by
$$ (a,\sigma)\cdot(a',\sigma')=(a\sigma(a'),\sigma\sigma'),\ \ \ a,a'\in G^n\ \ \text{and}\ \ \sigma,\sigma'\in S_n.$$
The groups $\tG_n$ form an inductive chain of group:
\begin{equation}\label{chain}
\{1\}\subset\tG_1\subset\tG_2\subset\dots\dots\subset\tG_{n-1}\subset\tG_n\subset\dots\ ,
\end{equation}
where the subgroup of $\tG_n$ isomorphic to $\tG_{n-1}$ is formed by the elements of the form $(a,\sigma)$, $a\in G^{n-1}$ and $\sigma\in S_{n-1}$. This allows to consider elements of $\C\tG_{n-1}$ as elements of $\C\tG_n$ and we will often do this without mentioning.

\vskip .2cm
We denote by $s_1,\dots,s_{n-1}$ the following elements of $\tG_n$:
\begin{equation}\label{si}
s_i=\Bigl(1_{G^n}\,,\,\pi_i\Bigr),\ \ i=1,\dots,n-1,
\end{equation}
where $1_{G^n}$ is the unit element of $G^n$ and $\pi_i$ is the transposition of $i$ and $i+1$. The elements $s_1,\dots,s_{n-1}$ satisfy the following relations:
\begin{equation}\label{rel-Sn}\begin{array}{ll}
s_i^2=1 & \text{for $i=1,\dots,n-1$,}\\
s_is_{i+1}s_i=s_{i+1}s_is_{i+1} & \text{for $i=1,\dots,n-2$,}\\
s_is_j=s_js_i & \text{for $i,j=1,\dots,n-1$ such that $|i-j|>1$,}
\end{array}\end{equation}
and generate a subgroup of $\tG_n$ isomorphic to the symmetric group $S_n$.

\vskip .2cm
 For $j=1,\dots,n$, let $\varpi_j$ be the injective morphism from $\C G$ to $\C G^n$ defined by 
\begin{equation}\label{inj-g}\varpi_j(g):=(1_G,\dots,1_G,\,g\,,1_G,\dots,1_G)\ \ \ \ \text{for $g\in \C G$,}
\end{equation}
where $1_G$ is the unit element of the group $G$ and, in the right hand side of (\ref{inj-g}), $g$ is in the $j$-th position. Let also $\iota$ be the natural injective morphism from $\C G^n$ to $\C\tG_n$, given by
$$\iota(a):=(a,1_{S_n})\ \ \ \ \text{for $a\in \C G^n$,}$$
where $1_{S_n}$ is the unit element of the symmetric group $S_n$. 

\vskip .2cm
For any $j=1,\dots,n$, the composition $\iota\circ\varpi_j$ is an injective morphism from $\C G$ to $\C \tG_n$. We will use the following notation:
\begin{equation}\label{inj-Gn}g_j:=\iota\circ\varpi_j(g)\ \ \ \text{for $g\in \C G$ and $j=1,\dots,n$.}\end{equation}
For $j=1,\dots,n$, we define, similarly to (\ref{bax-g}),
\begin{equation}\label{bax-gi}
g^{(\alpha)}_j(v):=\frac{\prod_{\xi^{(\alpha)}\in S^{(\alpha)}}(v-\xi^{(\alpha)})}{v-g_j^{(\alpha)}}\ \ \ \ \ \text{for $\alpha=1,\dots,m$\,.}
\end{equation}

\vskip .2cm
Now consider the following element of the group ring of $G\times G$:
\begin{equation}\label{def-e}\mathsf{e}:=\frac{1}{|G|}\sum_{g\in G}(g,g^{-1})\ .\end{equation}
As noticed in \cite{Pus}, the element $\mathsf{e}$ satisfies
\begin{equation}\label{rel-b}
(g,h)\,\mathsf{e}=\mathsf{e}\,(h,g)\ \ \ \ \ \text{for any $g,h\in G$\,,}
\end{equation}
and this implies that $\mathsf{e}$ acts as follows in the irreducible representations of $G\times G$:
\begin{equation}\label{eq-idem-e}\rho_{\nu}\otimes\rho_{\nu'}(\mathsf{e})=\frac{\delta_{\nu,\nu'}}{\di_{\nu}}\,\text{P}_{W_{\nu}\otimes W_{\nu}}\ \ \ \ \text{for $\nu,\nu'=1,\dots,m$\,,}\end{equation}
where $\text{P}_{W_{\nu}\otimes W_{\nu}}$ is the permutation operator of the space $W_{\nu}\otimes W_{\nu}$ ($\text{P}(u\otimes v)=v\otimes u$ for $u,v\in W_{\nu}$). 

\paragraph{Elements $e_{i,j}$.} For any $i,j=1,\dots,n$ such that $i\neq j$, the map $\varpi_i\otimes\varpi_j$ is an injective morphism from the group ring of $G\times G$ to $\C G^n$, and thus the composition $\iota\circ(\varpi_i\otimes\varpi_j)$ is an injective morphism from the group ring of $G\times G$ to $\C\tG_n$. We define, for $i=1,\dots,n$, $e_{i,i}:=1$ and
\begin{equation}\label{def-eij}
e_{i,j}:=\iota\bigl(\,\varpi_i\otimes\varpi_j(\mathsf{e})\,\bigr)\ \ \ \ \ \text{for $i,j=1,\dots,n$ such that $i\neq j$,}
\end{equation}
With the notation (\ref{inj-Gn}), elements $e_{i,j}$ can be written as
\[e_{i,j}=\frac{1}{|G|}\sum_{g\in G}g_i^{\phantom{-1}}\!\!\!g_j^{-1}\ ,\ \ \ \ i,j=1,\dots,n\ .\]
 By construction, we have $e_{i,j}=e_{j,i}$ for $i,j=1,\dots,n$, and
\begin{equation}\label{eq-s-e}s_k\,e_{i,j}=e_{\pi_k(i),\pi_k(j)}\,s_k\ \ \ \ \text{for $i,j=1,\dots,n$ and $k=1,\dots,n-1$\,,}\end{equation}
where we recall that $\pi_k$ is the transposition of $k$ and $k+1$. 

\vskip .2cm
For $1\leq i<j\leq n$, the element $\varpi_i\otimes\varpi_j(\mathsf{e})$ acts in the irreducible representations of $G^n$ as follows (this is implied by (\ref{eq-idem-e}))
\begin{equation}\label{eq-eij-rep}
\rho_{\nu_1}\otimes\dots\otimes\rho_{\nu_n}\bigl(\varpi_i\otimes\varpi_j(\mathsf{e})\bigr)=\frac{\delta_{\nu_i,\nu_j}}{\di_{\nu_i}}\,\text{P}_{i,j}\ ,\ \ \ \ \text{for any $\nu_1,\dots,\nu_n\in\{1,\dots,m\}$,}
\end{equation}
where $\text{P}_{i,j}$ is the operator on $W_{\nu_1}\otimes\dots\otimes W_{\nu_n}$ which permutes the $i$-th and $j$-th spaces and acts as the identity operator anywhere else, that is 
$$\text{P}_{i,j}(u_1\otimes \dots u_i \dots u_j \,.\,.\,.\, \otimes u_n)=u_1\otimes \dots u_j \dots u_i \,.\,.\,.\, \otimes u_n\ ,\qquad\text{$u_a\in W_{\nu_a}$ for $a=1,\dots,n\,.$}$$

The following relations are consequences of (\ref{def-eij}) and (\ref{eq-eij-rep}):
\begin{equation}\label{eq-e-e}
\begin{array}{l} e_{i,j}\,e_{k,l}\ =\ e_{k,l}\,e_{i,j}\quad\ \ \ \ \text{for $\,1\leq i<j<k<l\leq n$\, such that $k-j>1$\,,}\\[0.4em]
e_{i,i+1}e_{k,l}=e_{\pi_i(k),\pi_i(l)}e_{i,i+1}\quad\ \ \ \ \text{for $i=1,\dots,n-1$ and $k,l=1,\dots,n$\,,}\\[0.4em]
e_{i,i+1}\,e_{i+1,i+2}\,e_{i,i+1}\ =\ e_{i+1,i+2}\,e_{i,i+1}\,e_{i+1,i+2}\quad\ \ \ \ \text{for $i=1,\dots,n-2$\,.}\end{array}
\end{equation}

\section{{\hspace{-0.55cm}.\hspace{0.55cm}}Jucys--Murphy elements}\label{sec-JM}

The Jucys--Murphy elements $\tj_1,\dots,\tj_n$ for the group $\tG_n$ are defined \cite{Pus,Wan} by the following initial condition and recursion:
\begin{equation}\label{def-jm} \tj_1=0\ \ \ \ \text{and}\ \ \ \ \tj_{i+1}=s_i\tj_is_i+e_{i,i+1}s_i\ \ \ \text{for $i=1,\dots,n-1$\,.}\end{equation}
The Jucys--Murphy elements $\tj_1,\dots,\tj_n$ form a commutative set of elements,
\begin{equation}\label{comm-JM}\tj_k\,\tj_l=\tj_l\,\tj_k\ \ \ \ \ \text{for $k,l=1,\dots,n$\,;}\end{equation}
moreover they satisfy, for $k=1,\dots,n$,
\begin{equation}\label{eq-j-s-g}\begin{array}{l}s_i\,\tj_k=\tj_k\,s_i\ \ \ \ \ \text{for $i=1,\dots,n-1$ such that $i\neq k-1,k$\,,}\\[0.4em]
g_l\,\tj_k=\tj_k\,g_l\ \ \ \ \ \text{for $l=1,\dots,n$ and for any $g\in G$\,.}
\end{array}\end{equation}

Recall that, for $\alpha=1,\dots,m$, the element $g^{(\alpha)}$ is the central element of $\C G$ defined by (\ref{def-g}), and that, for $l=1,\dots,n$, the element $g_l^{(\alpha)}$ is the image in $\C \tG_n$ of $g^{(\alpha)}$ by the injective morphism $\iota\circ\varpi_l$.
The relation (\ref{comm-JM}) and the second relation in (\ref{eq-j-s-g}) imply that the set of elements $\{g^{(\alpha)}_l,\ \alpha=1,\dots,m,\ l=1,\dots,n\}$ together with the Jucys--Murphy elements $\tj_1,\dots,\tj_n$ form a commutative set of elements of $\C\tG_n$.

\vskip .2cm
We state here the following Lemma, which will be used later in the proof of the fusion formula in Section \ref{sec-fus}.
\begin{lemm}\label{lemm-tj}
For $l=1,\dots,n-1$, we have
\begin{equation}\label{eq-tj}
\tj_l\,s_ls_{l+1}\dots s_{n-1}=s_ls_{l+1}\dots s_{n-1}\,\tj_n-\sum_{k=l}^{n-1}s_ls_{l+1}\dots \widehat{s_k}\dots s_{n-1}\,e_{k,n}\ ,
\end{equation}
where $s_ls_{l+1}\dots \widehat{s_k}\dots s_{n-1}$ stands for the product $s_ls_{l+1}\dots s_{n-1}$ with $s_k$ removed.
\end{lemm}
\emph{Proof.} We prove the formula (\ref{eq-tj}) by induction on $n-l$. The basis of induction (for $l=n-1$) is
\[\tj_{n-1}s_{n-1}=s_{n-1}\tj_n-e_{n-1,n}\ ,\]
which comes immediately from the recurrence relation in (\ref{def-jm}).

\vskip .2cm
For $l<n-1$,we  use $\tj_l\,s_l=s_l\,\tj_{l+1}-e_{l,l+1}$ to write
\[\tj_l\,s_ls_{l+1}\dots s_{n-1}=s_l\,\tj_{l+1}\,s_{l+1}\dots s_{n-1}-e_{l,l+1}\,s_{l+1}\dots s_{n-1}\ .\]
The formula (\ref{eq-tj}) follows, using the induction hypothesis and $e_{l,l+1}\,s_{l+1}\dots s_{n-1}=s_{l+1}\dots s_{n-1}\,e_{l,n}$.\hfill$\square$

\section{{\hspace{-0.55cm}.\hspace{0.55cm}}Baxterized elements}\label{sec-bax}

Recall the standard Artin presentation of the symmetric group $S_n$: the group $S_n$ is generated by elements $\overline{s}_1,\dots\overline{s}_{n-1}$ with defining relations as in (\ref{rel-Sn}), with $s_i$ replaced by $\overline{s}_i$.
 The standard Baxterization for the generators $\overline{s}_1,\dots,\overline{s}_{n-1}$ of $S_n$ is:
\[ \overline{s}_i(c,c'):=\overline{s}_i+\frac{1}{c-c'}\ \ \ \ \ \text{for $i=1,\dots,n-1$\,;}\]
the parameters $c$ and $c'$ are called \emph{spectral} parameters.

\vskip .2cm
The elements $s_1,\dots,s_{n-1}$, defined by (\ref{si}), of $\tG_n$ generate a subgroup isomorphic to the symmetric group $S_n$. We define a generalization of the Baxterized elements $\overline{s}_i(c,c')$, associated to the group $\tG_n$, as follows:
\begin{equation}\label{def-baxt}
s_i(c,c'):=s_i+\frac{e_{i,i+1}}{c-c'}\ \ \ \ \ \text{for $i=1,\dots,n-1$\,.}
\end{equation}

\begin{prop}\label{prop-baxt}
The Baxterized elements $s_i(c,c')$ satisfy the Yang--Baxter equation with spectral parameters:
\begin{equation}\label{YB}
s_i(c,c')s_{i+1}(c,c'')s_i(c',c'')=s_{i+1}(c',c'')s_i(c,c'')s_{i+1}(c,c')\ \ \ \ \ \text{for $i=1,\dots,n-2$\,;}
\end{equation}
they satisfy also
\begin{equation}\label{YB2}
\begin{array}{l}s_i(c,c')s_j(d,d')=s_j(d,d')s_i(c,c')\ \ \ \ \ \text{for $i,j=1,\dots,n-1$ such that $|i-j|>1$\,,}\\[0.4em]
s_i(c,c')s_i(c',c)=1-\displaystyle\frac{e_{i,i+1}^2}{(c-c')^2}\ \ \ \ \ \text{for $i=1,\dots,n-1$\,.}
\end{array}
\end{equation}
\end{prop}
\emph{Proof.} The first relation in (\ref{YB2}) is immediate using that $s_is_j=s_js_i$ and $e_{i,i+1}\,e_{j,j+1}=e_{j,j+1}\,e_{i,i+1}$ for $i,j=1,\dots,n-1$ such that $|i-j|>1$ (see (\ref{rel-Sn}) and (\ref{eq-e-e})).

\vskip .1cm
The second relation in (\ref{YB2}) follows from a direct calculation in which one uses that $s_i\,e_{i,i+1}=e_{i,i+1}\,s_i$ for $i=1,\dots,n-1$ (see (\ref{eq-s-e})).

\vskip .1cm
To finish the proof of the Proposition, we develop successively both sides of (\ref{YB}) and, using (\ref{eq-s-e}), we move in each term all elements $e_{k,l}$ on the left of the elements $s_i$ and $s_{i+1}$. We compare the two results:

\vskip .1cm
$-$ The cubic terms in the elements $s_i$ and $s_{i+1}$ coincide since $s_is_{i+1}s_i=s_{i+1}s_is_{i+1}$.

\vskip .1cm
$-$ In both sides, the coefficient in front of $s_is_{i+1}$ is $\displaystyle\frac{e_{i+1,i+2}}{c'-c''}\,$.

\vskip .1cm
$-$ In both sides, the coefficient in front of $s_{i+1}s_i$ is $\displaystyle\frac{e_{i,i+1}}{c-c'}\,$.

\vskip .1cm
$-$ In the left hand side, the coefficient in front of $s_i$ is $$\frac{e_{i,i+1}\,e_{i+1,i+2}}{(c-c')(c-c'')}+\frac{e_{i,i+2}\,e_{i,i+1}}{(c-c'')(c'-c'')}\ .$$
Using the second line in (\ref{eq-e-e}), we transform this expression in $\displaystyle\frac{e_{i,i+1}\,e_{i+1,i+2}}{(c-c')(c'-c'')}\,$, which coincides with the coefficient in front of $s_i$ in the right hand side (we also use (\ref{eq-e-e}) here).

\vskip .1cm
$-$ A similar calculation shows that the coefficient in front of $s_{i+1}$ in the right hand side is equal to $\displaystyle\frac{e_{i+1,i+2}\,e_{i,i+1}}{(c-c')(c'-c'')}\,$, which coincides with the coefficient in front of $s_{i+1}$ in the left hand side.

\vskip .1cm
$-$ The remaining terms in the left hand side are $$\frac{e_{i,i+2}}{c-c''}+\frac{e_{i,i+1}\,e_{i+1,i+2}\,e_{i,i+1}}{(c-c')(c-c'')(c'-c'')}\ ,$$ which coincide with the remaining terms 
$$ \frac{e_{i,i+2}}{c-c''}+\frac{e_{i+1,i+2}\,e_{i,i+1}\,e_{i+1,i+2}}{(c'-c'')(c-c'')(c-c')}\ $$
of the right hand side, using the third line in (\ref{eq-e-e}).
\hfill$\square$

\vskip .2cm
\textbf{Remark.} For $i=1,\dots,n-1$ and $\alpha=1,\dots,m$, the elements $s_i(c,c')$ and $g_i^{(\alpha)}(v)$ satisfy a certain limit of the reflection equation with spectral parameters (see for example \cite{IO}), namely
\begin{equation}\label{refl}
s_i(c,c')\,g_i^{(\alpha)}\!(c)\,s_i\,g_i^{(\alpha)}\!(c')=g_i^{(\alpha)}\!(c')\,s_i\,g_i^{(\alpha)}\!(c)\,s_i(c,c')\ .
\end{equation}
Indeed, due to (\ref{bax-gi}) and (\ref{YB2}) and the fact that $e_{i,i+1}$ commutes with $s_i$ and $g_i^{(\alpha)}$, the equality (\ref{refl}) is equivalent to
\[s_i(c',c)\,(c-g_i^{(\alpha)})\,s_i\,(c'-g_i^{(\alpha)})=(c'-g_i^{(\alpha)})\,s_i\,(c-g_i^{(\alpha)})\,s_i(c',c)\ ,\]
which is proved by a straightforward calculation. We skip the details and just indicate that one uses that $g_i^{(\alpha)}e_{i,i+1}=g_{i+1}^{(\alpha)}e_{i,i+1}$, which follows from (\ref{rel-b}) together with the fact that $g_i^{(\alpha)}$ commutes with $e_{i,i+1}$.

\section{{\hspace{-0.55cm}.\hspace{0.55cm}}Standard $m$-tableaux and idempotents of $\tG_n$} \label{sec-tab}

\subsection{{\hspace{-0.50cm}.\hspace{0.50cm}}$m$-partitions and $m$-tableaux}

\paragraph{1. $m$-partitions.}
Let $\lambda\vdash n$ be a partition of $n$ (we shall also say \emph{of size $n$}), that is, $\lambda=(\lambda_1,\dots,\lambda_k)$, where $\lambda _j$, $j=1,\dots,k$, are positive integers, $\lambda_1\geq\lambda_2\geq\dots\geq\lambda_k$ and the size of $\lambda$ is $|\lambda|:=\lambda_1+\dots+\lambda_k=n$. 

\vskip .2cm
We identify partitions with their Young diagrams: the Young diagram of $\lambda$ is a left-justified array of rows of nodes containing $\lambda_j$ nodes in the $j$-th row, $j=1,\dots,k$; the rows are numbered from top to bottom. Note that the number of nodes in the diagram of $\lambda$ is equal to $n$, the size of $\lambda$.

Recall that, for a partition $\lambda$, a node $\theta$ of $\lambda$ is called {\it removable} if the set of nodes obtained from $\lambda$ by removing $\theta$ is still a partition. A node $\theta'$ not in $\lambda$ is called {\it  addable} if the set of nodes obtained from $\lambda$ by adding $\theta'$ is still a partition.

\vskip .2cm
An $m$-partition $\lam$, or a Young $m$-diagram $\lam$, of size $n$ is an $m$-tuple of partitions, $$\lam:=(\lambda^{(1)},\dots,\lambda^{(m)})\,,$$ such that the total number of nodes in the associated Young diagrams is equal to $n$, that is $|\lambda^{(1)}|+\dots+|\lambda^{(m)}|=n$.

A pair $(\theta,k)$ consisting of a node $\theta$ and an integer $k\in\{1,\dots,m\}$ is called an $m$-node. The integer $k$ is called the position of the $m$-node.

\vskip .2cm
Let $\lam=(\lambda^{(1)},\dots,\lambda^{(m)})$ be an $m$-partition. An $m$-node $\thm=(\theta,k)\in\lam$ is called removable from $\lam$ if the node $\theta$ is removable from $\lambda^{(k)}$. An $m$-node $\thm'=(\theta',k')\notin\lam$ is called addable to $\lam$ if the node $\theta'$ is addable to $\lambda^{(k')}$. The set of $m$-nodes removable from $\lam$ is denoted by ${\cal{E}}_-(\lam)$ and the set of $m$-nodes addable to $\lam$ is denoted by ${\cal{E}}_+(\lam)$. For example, the removable/addable $3$-nodes (marked with -/+) for the $3$-partition $\left(\Box\!\Box,\varnothing,\Box\right)$
are
\[\left(\begin{array}{l}\fbox{$\phantom{-}$}\fbox{$-$}\fbox{$+$}\\ \fbox{$+$}\end{array}
\, ,\, \begin{array}{l}\fbox{$+$}\\ \phantom{\fbox{$-$}}\end{array}\, ,\,\begin{array}{l}\fbox{$-$}\fbox{$+$}\\ \fbox{$+$}\end{array}\right) \]

\vskip .2cm
For an $m$-node $\thm$ lying in the line $x$ and the column $y$ of the $k$-th diagram, we define $\pos(\thm):=k$ and $\cc(\thm):=y-x$. The number $\pos(\thm)$ is the position of $\thm$ and the number $\cc(\thm)$ is called the classical content of the $m$-node $\thm$. 

\vskip .2cm
For an $m$-partition $\lam$, we define
\begin{equation}\label{def-fg}
\mathsf{F}^G_{\lam}:=\prod_{\thm\in\lam}
\Biggl(\ \prod_{\alpha=1}^{m}\Bigl(\!\!\prod_{\textrm{\scriptsize{$\begin{array}{c}\xi^{(\alpha)}\in S^{(\alpha)}\\[0em]
\xi^{(\alpha)}\neq\xi_{\pos(\thm)}^{(\alpha)}\end{array}$}}
}\!\!\!\!\bigl(\xi_{\pos(\thm)}^{(\alpha)}-\xi^{(\alpha)}\bigr)\ \Bigr)\,\Biggr)
\ ,
\end{equation}
where the sets $S^{(\alpha)}$ and the numbers  $\xi_{\nu}^{(\alpha)}$, for $\alpha,\nu=1,\dots,m$, are defined in Section \ref{sec-def}.

\paragraph{2. Hook length.} 
Let $\lambda$ be a partition and $\theta$ be a node of $\lambda$. The hook of $\theta$ in $\lambda$ is the set of nodes of  $\lambda$ consisting of the node $\theta$ together with the nodes which lie either under $\theta$ in the same column or to the right of $\theta$ in the same row; the hook length $h_{\lambda}(\theta)$ is the cardinality of the hook of $\theta$ in $\lambda$. 

\vskip .2cm
The definition of the hook length is extended to $m$-partitions and $m$-nodes as follows. Let $\lam=(\lambda^{(1)},\dots,\lambda^{(m)})$ be an $m$-partition and $\thm=(\theta,k)$ an $m$-node of $\lam$. The hook length $h_{\lam}(\thm)$ of $\thm$ in $\lam$ is the hook length of the node $\theta$ in the $k$-th partition of $\lam$, that is
$$h_{\lam}(\thm):=h_{\lambda^{(k)}}(\theta)\ .$$

\vskip .2cm
For an $m$-partition $\lam=(\lambda^{(1)},\dots,\lambda^{(m)})$, we define
\begin{equation}\label{def-f}
\mathsf{F}_{\lam}:=\prod_{\thm\in\lam}h_{\lam}(\thm)=\prod_{k=1}^m\prod_{\theta\in\lambda^{(m)}}h_{\lambda^{(m)}}(\theta)\ .
\end{equation}

\paragraph{3. Standard $m$-tableaux.}
Let $\lam$ be an $m$-partition of $n$. An $m$-tableau of shape $\lam$ is a bijection between the set $\{1,\dots,n\}$ and the set of $m$-nodes in $\lam$, that is, an $m$-tableau of shape $\lam$ is obtained by placing the numbers $1,\dots,n$ in the $m$-nodes of $\lam$. We call the number $n$ the \emph{size} of the $m$-tableau. An $m$-tableau is standard if the numbers increase along any row and any column of every diagram in $\lam$. 

\vskip .2cm
For a standard $m$-tableau ${\cal{T}}$, we denote respectively by $\cc({\cal{T}}|i)$ and $\pos({\cal{T}}|i)$ the classical content and the position of the $m$-node with number $i$. For example, for the standard $3$-tableau ${\cal{T}}^{^{\phantom{A}}}\!\!\!\!={\textrm{$\left(
\,\fbox{\scriptsize{$1$}}\fbox{\scriptsize{$3$}}\, ,\,\varnothing\, ,\,\fbox{\scriptsize{$2$}}\,\right)$}}$, we have
\[\cc({\cal{T}}|1)=0\,,\ \ \cc({\cal{T}}|2)=0\,,\ \ \cc({\cal{T}}|3)=1\ \ \ \ \ \text{and}\ \ \ \ \  \pos({\cal{T}}|1)=1\,,\ \ \pos({\cal{T}}|2)=3\,,\ \ \pos({\cal{T}}|3)=1\,.\]

For a standard $m$-tableau ${\cal{T}}$, we define the \emph{$G$-content} $\cc^G({\cal{T}}|i)$ of the $m$-node with number $i$ by
\begin{equation}\label{G-cont}
\cc^G({\cal{T}}|i):=\frac{\cc({\cal{T}}|i)}{\di_{\pos({\cal{T}}|i)}}\ ,
\end{equation}
where we recall that $\di_1,\dots,\di_m$ are the dimensions of the pairwise non-isomorphic irreducible representations $W_1,\dots,W_m$ of $G$. For example, for $m=3$ and for the same standard $3$-tableau as above, we have:
\[\cc^G({\cal{T}}|1)=0\,,\ \ \cc^G({\cal{T}}|2)=0\,,\ \ \cc^G({\cal{T}}|3)=\frac{1}{\di_1}\ .\] 
 
\vskip .2cm
Let $\vm$ be an $m$-tuple of variables, $\vm:=(v^{(1)},\dots,v^{(m)})$. Let $N$ be a non-negative integer, $\lam$ an $m$-partition of size $N$ and ${\cal{T}}$  a standard $m$-tableau  of shape $\lam$. For brevity, set $\cc_i:=\cc({\cal{T}}|i)$ and $\pos_i:=\pos({\cal{T}}|i)$ for $i=1,\dots,N$. We define
\begin{equation}\label{def-F1-G}
F^G_{_{{\cal{T}}}}(\vm):=\Bigl(\prod_{\alpha=1}^{m}\!\!\!\prod_{\textrm{\scriptsize{$\begin{array}{c}\xi^{(\alpha)}\in S^{(\alpha)}\\[0em]
\xi^{(\alpha)}\neq\xi_{\pos_N}^{(\alpha)}\end{array}$}}
}\frac{1}{v^{(\alpha)}-\xi^{(\alpha)}}\Bigr)\ ,
\end{equation}
and
\begin{equation}\label{def-F1}
F_{_{{\cal{T}}}}(u):=\frac{u-\cc_N}{u}\,\prod_{i=1}^{N-1}\frac{(u-\cc_i)^2}{(u-\cc_i)^2-\delta_{\pos_i,\pos_N}}\ ,
\end{equation}
where $\delta_{\pos_i,\pos_N}$ is the Kronecker delta.

\vskip .2cm
Let $\mum$ be the shape of the standard $m$-tableau obtained from ${\cal{T}}$ by removing the $m$-node containing the number $N$.
Then $F^G_{_{{\cal{T}}}}(\vm)$ is non-singular at $v^{(\alpha)}=\xi_{\pos_N}^{(\alpha)}$, $\alpha=1,\dots,m$, and moreover, from (\ref{def-fg}),
\begin{equation}\label{eq-F1-G}
F^G_{_{{\cal{T}}}}(\vm)\Bigr\rvert_{v^{(\alpha)}=\xi_{\pos_N}^{(\alpha)}\,,\ \alpha=1,\dots,m}=\left(\mathsf{F}^G_{\lam}\right)^{-1}\mathsf{F}_{\mum}^{G}\ .
\end{equation}

\vskip .2cm
We will also need the following known result \cite{Mo,OPdA} concerning the function $F_{_{{\cal{T}}}}(u)$.
\begin{lemm}
\label{lemm-f}
The rational function $F_{_{{\cal{T}}}}(u)$ is non-singular at $u=\cc_{N}$ and moreover
\begin{equation}\label{eq-F1}F_{_{{\cal{T}}}}(u)\Bigr\rvert_{u=\cc_N}=\mathsf{F}_{\lam}^{-1}\mathsf{F}^{\phantom{-1}}_{\mum}\ .
\end{equation}
\end{lemm}

\subsection{{\hspace{-0.50cm}.\hspace{0.50cm}}Idempotents of $\tG_n$ corresponding to standard $m$-tableaux}
 
 The irreducible representations of $\tG_n$ are indexed by the set of $m$-partitions of $n$ (see \emph{e.g.} \cite{JK} or \cite{Mac}). Let $\lam$ be an $m$-partition of $n$ and denote by $V_{\lam}$ the vector space carrying the irreducible representation of $\tG_n$ corresponding to $\lam$. Let ${\cal{T}}$ be a standard $m$-tableau of shape $\lam$. We set
 \[W_{{\cal{T}}}:=\bigotimes_{i=1}^nW_{\pos({\cal{T}}|i)}\ ,\]
 where we recall that $W_1,\dots,W_m$ are the irreducible representation spaces of the group $G$.

 Then the vector space $V_{\lam}$ admits the following decomposition:
 \[V_{\lam}=\bigoplus W_{{\cal{T}}}\ ,\]
 where the direct sum is over the set of standard $m$-tableaux of shape $\lam$. We denote by $E_{\cal{T}}$ the idempotent of $\C\tG_n$ corresponding to the subspace $W_{{\cal{T}}}$.
 
\vskip .2cm
Recall that the set of elements $\{g^{(\alpha)}_k,\ \alpha=1,\dots,m,\ k=1,\dots,n\}$ together with the Jucys--Murphy elements $\tj_1,\dots,\tj_n$ form a commutative family of elements of $\C\tG_n$ (see Section \ref{sec-JM}). For any standard $m$-tableau  ${\cal{T}}$ of size $n$, the subspace $W_{{\cal{T}}}$ is a common eigenspace for this family of elements. Moreover, we have:
\begin{equation}\label{spec-g}
g^{(\alpha)}_k E_{{\cal{T}}}=E_{{\cal{T}}}\,g^{(\alpha)}_k=\xi_{\pos({\cal{T}}|k)}^{(\alpha)}E_{{\cal{T}}}\ \ \ \ \text{for $\alpha=1,\dots,m$ and $k=1,\dots,n$,}
\end{equation}
where the numbers $\xi_i^{\alpha)}$ are defined by formula (\ref{def-p}); from the results in \cite{Pus}, we also have
\begin{equation}\label{spec-j}
\tj_k E_{{\cal{T}}}=E_{{\cal{T}}} \tj_k=\cc^G({\cal{T}}|k)E_{{\cal{T}}}\ \ \ \ \text{for $k=1,\dots,n$.}
\end{equation}
For two standard $m$-tableaux ${\cal{T}}$ and ${\cal{T}}'$ of size $n$ such that ${\cal{T}}\neq {\cal{T}}'$, we have:\\[0.1cm]
$-$ either there exists $\alpha\in\{1,\dots,m\}$ and $k\in\{1,\dots,n\}$ such that $\xi_{\pos({\cal{T}}|k)}^{(\alpha)}\neq\xi_{\pos({\cal{T}}'|k)}^{(\alpha)}$.\\[0.1cm]
$-$ or there exists $k\in\{1,\dots,n\}$ such that $\cc^G({\cal{T}}|k)\neq\cc^G({\cal{T}}'|k)$.\\[0.1cm]
Indeed, if there is some $k\in\{1,\dots,n\}$ such that $\pos({\cal{T}}|k)\neq\pos({\cal{T}}'|k)$ then $\xi_{\pos({\cal{T}}|k)}^{(\alpha)}\neq\xi_{\pos({\cal{T}}'|k)}^{(\alpha)}$ for some $\alpha\in\{1,\dots,m\}$ (see Section \ref{sec-def}).  If $\pos({\cal{T}}|k)=\pos({\cal{T}}'|k)$ for any $k=1,\dots,m$ then, unless $\cc^G({\cal{T}}|k)\neq\cc^G({\cal{T}}'|k)$ for some $k\in\{1,\dots,n\}$,
we have that ${\mathcal{T}}$ and ${\mathcal{T}}'$ must have the same shape; moreover the $m$-tableau ${\mathcal{T}}'$ must be obtained from ${\mathcal{T}}$ by permuting the entries inside each diagonal of each diagram. As both $m$-tableaux are standard, we must have ${\cal{T}}={\cal{T}}'$ which contradicts the assumption. Thus $\cc^G({\cal{T}}|k)\neq\cc^G({\cal{T}}'|k)$ for some $k\in\{1,\dots,n\}$.

Thus, the idempotent $E_{{\cal{T}}}$ can be expressed in terms of the elements $g^{(\alpha)}_k,\ \alpha=1,\dots,m,\ k=1,\dots,n$ and the Jucys--Murphy elements $\tj_1,\dots,\tj_n$.

\vskip .2cm
Let $\lam$ be an $m$-partition of $n$ and ${\cal{T}}$ be a standard $m$-tableau of shape $\lam$. Denote by $\thm$ the $m$-node of ${\cal{T}}$ containing the number $n$, and for brevity, set $\cc^G_n:=\cc^G(\thm)$ and $\pos_n=\pos(\thm)$. As the $m$-tableau ${\cal{T}}$ is standard, the $m$-node $\thm$ of $\lam$ is removable. Let ${\cal{U}}$ be the standard $m$-tableau obtained from ${\cal{T}}$ by removing $\thm$ and let $\mum$ be the shape of ${\cal{U}}$.

We have the following inductive formula for $E_{{\cal{T}}}$:
\begin{equation}\label{idem-JM}
E_{{\cal{T}}}=E_{{\cal{U}}}\!\!\!\prod_{\kam\colon\!\!\!\! \begin{array}{l}\scriptstyle{\kam\,\in\, 
{\cal{E}}_+(\mum)}\\\scriptstyle{\cc^G(\kam)\neq \cc^G_n}\end{array}  }\hspace{-0.4cm} \frac{\tj_{n}-\cc^G(\kam)}{\cc^G_n-\cc^G(\kam)}
\ \ \prod_{\alpha=1}^m\Biggl(\ \ \prod_{\kam\colon\!\!\!\! \begin{array}{l}\scriptstyle{\kam\,\in\, 
{\cal{E}}_+(\mum)}\\\scriptstyle{\xi^{(\alpha)}_{\pos(\kam)}\neq \xi^{(\alpha)}_{\pos_n}}\end{array}  }\hspace{-0.3cm}\frac{g^{(\alpha)}_n-\xi_{\pos(\kam)}^{(\alpha)}}{\xi_{\pos_n}^{(\alpha)}-\xi_{\pos(\kam)}^{(\alpha)}}\Biggr)\ ,
\end{equation}
where we recall that ${\cal{E}}_+(\mum)$ is the set of addable $m$-nodes of the $m$-partition $\mum$; the initial condition is $E_{{\cal{U}}_0}=1$, where ${\cal{U}}_0$ is the unique standard $m$-tableau of size $0$. Note that the element $E_{{\cal{U}}}$ in (\ref{idem-JM}) is an idempotent of $\C\tG_{n-1}$ and we consider it as an element of $\C\tG_n$ due to the chain property of the groups $\tG_n$ (see (\ref{chain})).

\vskip .2cm
Let $\{ {\cal{T}}_1,\dots,{\cal{T}}_k\}$ be the set of pairwise different standard $m$-tableaux that can be obtained from ${\cal{U}}$ by adding an $m$-node containing the number $n$. Note that ${\cal{T}}\in\{ {\cal{T}}_1,\dots,{\cal{T}}_k\}$. We have:
\begin{equation}\label{Eu-Et}E_{{\cal{U}}}=\sum_{i=1}^k E_{{\cal{T}}_i}\ .\end{equation}
Consider the following rational function in $u$ and $\vm$ with values in $\C \tG_n$
\begin{equation}\label{Eu-Et2} E_{{\cal{U}}}\ \ \frac{u-\cc^G_n}{u-\tj_n}\ \ \prod_{\alpha=1}^m\ \frac{v^{(\alpha)}-\xi_{\pos_n}^{(\alpha)}}{v^{(\alpha)}-g^{(\alpha)}_n}\ ,\end{equation}
and replace $E_{{\cal{U}}}$ by the right hand side of (\ref{Eu-Et}). Then formulas (\ref{spec-g})--(\ref{spec-j}) imply that the rational function (\ref{Eu-Et2}) is non-singular at $u=\cc^G_n$ and $v^{(\alpha)}=\xi_{\pos_n}^{(\alpha)}$, $\alpha=1,\dots,m$, and moreover, 
\begin{equation}\label{idem-JM-fin} E_{{\cal{U}}}\ \ \frac{u-\cc^G_n}{u-\tj_n}\ \ \prod_{\alpha=1}^m\ \frac{v^{(\alpha)}-\xi_{\pos_n}^{(\alpha)}}{v^{(\alpha)}-g^{(\alpha)}_n}\ 
\Bigr\rvert_{\hspace{-0.15cm}\textrm{\scriptsize{$\begin{array}{l}u=\cc^G_n\\[-0.15em]v^{(\alpha)}=\xi_{\pos_n}^{(\alpha)}\,,\ \alpha=1,\dots,m\end{array}$}}}=E_{{\cal{T}}}\ .\end{equation}

 \vskip .2cm
 \textbf{Remarks.} \textbf{(i)} The elements $E_{\cal{T}}$, with ${\cal{T}}$ running through the set of standard $m$-tableaux of size $n$, form a complete system of pairwise orthogonal idempotents of $\C \tG_n$. For a standard $m$-tableau ${\cal{T}}$, the idempotent $E_{\cal{T}}$ is primitive if and only if $W_{\pos({\cal{T}}|i)}$ is one-dimensional for any $i\in\{1,\dots,n\}$. Thus the elements $E_{\cal{T}}$ form a complete system of primitive pairwise orthogonal idempotents of $\C \tG_n$ if and only if the group $G$ is Abelian.
 
\vskip .2cm
\textbf{(ii)} The preceding remark can also be expressed as follows: The commutative subalgebra generated by the set $\{g^{(\alpha)}_k,\ \alpha=1,\dots,m,\ k=1,\dots,n\}$ together with the Jucys--Murphy elements $\tj_1,\dots,\tj_n$ is a maximal commutative subalgebra of $\C \tG_n$ if and only if the group $G$ is Abelian.\hfill$\square$

\section{{\hspace{-0.55cm}.\hspace{0.55cm}}Fusion formula}\label{sec-fus}

Let 
\begin{equation}\label{def-gamma}\Gamma(\underline{v_1},\dots,\underline{v_n}):=\displaystyle\prod_{i=1}^n\prod_{\alpha=1}^mg^{(\alpha)}_i(v^{(\alpha)}_i)\ ,\end{equation}
where the functions $g^{(\alpha)}_i(v_i^{(\alpha)})$, $\alpha=1,\dots,m$ and $i=1,\dots,n$, are defined by (\ref{bax-gi}). 

\vskip .2cm
Set $\phi_{1}(u):=1$ and, for $k=2,\dots,n$, let
\begin{equation}\label{def-phi1}\begin{array}{l}\phi_{k}(u_1,\dots,u_{k-1},u):=
s_{k-1}(u,u_{k-1})\phi_{k-1}(u_1,\dots,u_{k-2},u)s_{k-1}\\[.2cm]
\hspace{2cm}=s_{k-1}(u,u_{k-1})s_{k-2}(u,u_{k-2})\dots s_1(u,u_1)\cdot s_1\dots s_{k-2}s_{k-1}\ .
\end{array}\end{equation}
We define the following rational function with values in $\C\tG_n$:
\begin{equation}\label{def-Phi1}
\Phi(u_1,\dots,u_n,\underline{v_1},\dots,\underline{v_n}):=\phi_n(u_1,\dots,u_n)\,\phi_{n-1}(u_1,\dots,u_{n-1})\dots\dots\phi_1(u_1)\,\Gamma(\underline{v_1},\dots,\underline{v_n})\ .
\end{equation}

\vskip .2cm
Let $\lam$ be an $m$-partition of size $n$ and ${\cal{T}}$ a standard $m$-tableau of shape $\lam$. For $i=1,\dots,n$, we set $\cc_i:=\cc({\cal{T}}|i)$, $\cc^G_i:=\cc^G({\cal{T}}|i)$ and $\pos_i:=\pos({\cal{T}}|i)$. Let also ${\cal{U}}$ be the standard $m$-tableau obtained from ${\cal{T}}$ by removing the $m$-node with number $n$, and let $\mum$ be the shape of ${\cal{U}}$. 

\begin{theo}\label{prop-fus1}
The idempotent $E_{{\cal{T}}}$ of $\C\tG_n$ corresponding to the standard $m$-tableau ${\cal{T}}$ can be obtained by the following consecutive evaluations
\begin{equation}\label{eq-idem-fin1}
E_{{\cal{T}}}=\frac{1}{\mathsf{F}^G_{\lam}\,\mathsf{F}^{\phantom{G}}_{\lam}}\,\Phi(u_1,\dots,u_n,\underline{v_1},\dots,\underline{v_n})\Bigr\rvert_{v^{(\alpha)}_i=\xi_{\pos_i}^{(\alpha)},\!\textrm{\scriptsize{$\begin{array}{l}i=1,\dots,n\\[0.1em]\alpha=1,\dots,m\end{array}$}}}\ \Bigr\rvert_{u_1=\cc^G_1}\dots\ \dots\ \Bigr\rvert_{u_n=\cc^G_n}\ .
\end{equation}
\end{theo}
\emph{Proof.} Define, for $k=1,\dots,n$,
\begin{equation}\label{def-phit}\widetilde{\phi}_k(u_1,\dots,u_{k-1},u,\vm):=\phi_{k}(u_1,\dots,u_{k-1},u)\prod_{\alpha=1}^mg^{(\alpha)}_k(v^{(\alpha)})\end{equation}
As $g_k^{(\alpha)}$ commutes with $s_i$ if $i<k-1$, we can rewrite the rational function $\Phi(u_1,\dots,u_n,\underline{v_1},\dots,\underline{v_n})$ as
\begin{equation}\label{Phi-ind}\Phi(u_1,\dots,u_n,\underline{v_1},\dots,\underline{v_n})=\widetilde{\phi}_n(u_1,\dots,u_n,\underline{v_n})\,\widetilde{\phi}_{n-1}(u_1,\dots,u_{n-1},\underline{v_{n-1}})\dots\dots\widetilde{\phi}_1(u_1,\underline{v_1})\ .\end{equation}

We prove the Theorem by induction on $n$. For $n=0$, there is nothing to prove.\\
For $n>0$, we use (\ref{Phi-ind}) and the induction hypothesis to rewrite the right hand side of (\ref{eq-idem-fin1}) as
\[\Bigl(\mathsf{F}^G_{\lam}\,\mathsf{F}^{\phantom{G}}_{\lam}\Bigr)^{-1}\mathsf{F}^G_{\mum}\,\mathsf{F}^{\phantom{G}}_{\mum}\cdot\widetilde{\phi}_n(\cc^G_1,\dots,\cc^G_{n-1},u_n,\underline{v_n})\,E_{{\cal{U}}}\Bigr\vert_{\hspace{-0.15cm}\textrm{\scriptsize{$\begin{array}{l}v_n^{(\alpha)}=\xi_{\pos_n}^{(\alpha)},\\
\alpha=1,\dots,m\end{array}$}}}\Bigr\rvert_{u_n=\cc^G_n}\ .\]
Now we use the Proposition \ref{prop-idem11} below to transform this expression into
\[\Bigl(\mathsf{F}^G_{\lam}\,\mathsf{F}^{\phantom{G}}_{\lam}\Bigr)^{-1}\mathsf{F}^G_{\mum}\,\mathsf{F}^{\phantom{G}}_{\mum}\cdot\Bigl(F^G_{_{{\cal{T}}}}(\vm)F_{_{{\cal{T}}}}(\di_{\pos_n}u)\Bigr)^{-1}\frac{u-\cc^G_n}{u-\tj_n}\ \ \prod_{\alpha=1}^m\ \frac{v^{(\alpha)}-\xi_{\pos_n}^{(\alpha)}}{v^{(\alpha)}-g^{(\alpha)}_n}\ E_{{\cal{U}}}\Bigr\vert_{\hspace{-0.15cm}\textrm{\scriptsize{$\begin{array}{l}v^{(\alpha)}=\xi_{\pos_n}^{(\alpha)},\\
\alpha=1,\dots,m\end{array}$}}}\Bigr\rvert_{u_n=\cc^G_n}\ .\]
We recall that $\di_{\pos_n}\cc^G_n=\cc_n$ and we use formulas (\ref{eq-F1-G}) and  (\ref{eq-F1}) concerning the functions $F^G_{_{{\cal{T}}}}(\vm)$ and $F_{_{{\cal{T}}}}(u)$, together with formula (\ref{idem-JM-fin}) to conclude.\hfill$\square$

\begin{prop}\label{prop-idem11}
Assume that $n\geq1$. We have 
\begin{equation}\label{eq-idem11}
F^G_{_{{\cal{T}}}}(\vm)F_{_{{\cal{T}}}}(\di_{\pos_n}u)\,\widetilde{\phi}_{n}(\cc^G_1,\dots,\cc^G_{n-1},u,\vm)\,E_{{\cal{U}}}\Bigr\vert_{\hspace{-0.15cm}\textrm{\scriptsize{$\begin{array}{l}v^{(\alpha)}=\xi_{\pos_n}^{(\alpha)},\\
\alpha=1,\dots,m\end{array}$}}}
\!\!\!\!=\,\frac{u-\cc^G_n}{u-\tj_n}\ \ \prod_{\alpha=1}^m\ \frac{v^{(\alpha)}-\xi_{\pos_n}^{(\alpha)}}{v^{(\alpha)}-g^{(\alpha)}_n}\ E_{{\cal{U}}}\Bigr\vert_{\hspace{-0.15cm}\textrm{\scriptsize{$\begin{array}{l}v^{(\alpha)}=\xi_{\pos_n}^{(\alpha)},\\
\alpha=1,\dots,m\end{array}$}}}\!\!.
\end{equation}
\end{prop}
\emph{Proof.} Notice that, from (\ref{bax-gi}) and (\ref{def-F1-G}),
\begin{equation}\label{def-FGfus}F^G_{_{{\cal{T}}}}(\vm)\prod_{\alpha=1}^mg^{(\alpha)}_n(v^{(\alpha)})=\prod_{\alpha=1}^m\ \frac{v^{(\alpha)}-\xi_{\pos_n}^{(\alpha)}}{v^{(\alpha)}-g^{(\alpha)}_n}\ .\end{equation}
Define
\begin{equation}\label{def-EUpn}E_{{\cal{U}},\pos_n}:=\prod_{\alpha=1}^m\ \frac{v^{(\alpha)}-\xi_{\pos_n}^{(\alpha)}}{v^{(\alpha)}-g^{(\alpha)}_n}E_{{\cal{U}}}\Bigr\vert_{\hspace{-0.15cm}\textrm{\scriptsize{$\begin{array}{l}v^{(\alpha)}=\xi_{\pos_n}^{(\alpha)},\\
\alpha=1,\dots,m\end{array}$}}}\!\!;
\end{equation}
the element $E_{{\cal{U}},\pos_n}$ is an idempotent which is equal to the sum of the idempotents $E_{{\cal{V}}}$, where ${\cal{V}}$ runs through the set of standard $m$-tableaux obtained from ${\cal{U}}$ by adding an $m$-node $\thm$ with number $n$ such that $\pos(\thm)=\pos_n$.

\vskip .2cm
Now, using (\ref{def-phit}), (\ref{def-FGfus}) and (\ref{def-EUpn}), the Proposition is a direct consequence of the Lemma \ref{lemm-idem11} below.\hfill$\square$

\begin{lemm}\label{lemm-idem11}
Assume that $n\geq1$. We have 
\begin{equation}\label{eq-idem12}
F_{_{{\cal{T}}}}(\di_{\pos_n}u)\,\phi_{n}(\cc^G_1,\dots,\cc^G_{n-1},u)\,E_{{\cal{U}},\pos_n}
=\frac{u-\cc^G_n}{u-\tj_n}E_{{\cal{U}},\pos_n}\ .
\end{equation}
\end{lemm}
\emph{Proof.} The left hand side of (\ref{eq-idem12}) is
\[F_{_{{\cal{T}}}}(\di_{\pos_n}u)(s_{n-1}+\frac{e_{n-1,n}}{u,\cc^G_{n-1}})\dots (s_1+\frac{e_{1,2}}{u-\cc^G_1})\cdot s_1\dots s_{n-1}\,E_{{\cal{U}},\pos_n}\ .\]
For $k=1,\dots,n-1$, one can verify, using relations (\ref{eq-s-e}) and (\ref{eq-e-e}), that 
$$\begin{array}{lc} & \displaystyle e_{k,k+1}\cdot (s_{k-1}+\frac{e_{k-1,k}}{u-\cc^G_{k-2}})\dots (s_1+\frac{e_{1,2}}{u-\cc^G_1})=(s_{k-1}+\frac{e_{k-1,k}}{u-\cc^G_{k-2}})\dots (s_1+\frac{e_{1,2}}{u-\cc^G_1})\cdot e_{1,k+1}\ ,\\[2em]
\text{and} & \displaystyle e_{1,k+1}\cdot s_1\dots s_{n-1}=s_1\dots s_{n-1}\cdot e_{k,n}\ .\end{array}$$
As moreover $e_{k,n}\,E_{{\cal{U}},\pos_n}=0$ if $\pos_k\neq\pos_n$ (see (\ref{eq-eij-rep})), the left hand side of (\ref{eq-idem12}) is equal to
\begin{equation}\label{lhs1}F_{_{{\cal{T}}}}(\di_{\pos_n}u)(s_{n-1}+\frac{\delta_{\pos_{n-1},\pos_n}e_{n-1,n}}{u-\cc^G_{n-1}})\dots (s_1+\frac{\delta_{\pos_1,\pos_n}e_{1,2}}{u-\cc^G_1})\cdot s_1\dots s_{n-1}\,E_{{\cal{U}},\pos_n}\ .\end{equation}

\vskip .2cm
 We prove the Lemma by induction on $n$. First assume that $\pos_i\neq\pos_n$ for $i=1,\dots,n-1$. In this situation, $E_{{\cal{U}},\pos_n}=E_{{\cal{T}}}$ and we have  $\cc^G_n=0$, $\tj_nE_{{\cal{T}}}=0$, $F_{_{{\cal{T}}}}(\di_{\pos_n}u)=1$ and, due to the formula (\ref{lhs1}), $F_{_{{\cal{T}}}}(\di_{\pos_n}u)\phi_n(\cc^G_1,\dots,\cc^G_{n-1},u)\,E_{{\cal{U}},\pos_n}=E_{{\cal{U}},\pos_n}$. So in this situation, the formula (\ref{eq-idem12}) is trivially satisfied; notice that, in particular, the basis of induction (for $n=1$) has been proved.
  
\vskip .2cm
Let $n>1$ and assume that there exists $l\in\{1,\dots,n-1\}$ such that $\pos_l=\pos_n$. Fix $l$ such that $\pos_l=\pos_n$ and $\pos_i\neq\pos_n$ for $i=l+1,\dots,n-1$. 

\vskip .2cm
Let ${\cal{V}}$ be the standard $m$-tableau obtained from ${\cal{U}}$ by removing the $m$-nodes with numbers $l+1,\dots,n-1$ and ${\cal{W}}$ be the standard $m$-tableau obtained from ${\cal{V}}$ by removing the $m$-node with number $l$. Define, similarly to (\ref{def-EUpn}),
\[E_{{\cal{W}},\pos_l}:=\prod_{\alpha=1}^m\ \frac{v^{(\alpha)}-\xi_{\pos_l}^{(\alpha)}}{v^{(\alpha)}-g^{(\alpha)}_{l}}E_{{\cal{W}}} \Bigr\vert_{\hspace{-0.15cm}\textrm{\scriptsize{$\begin{array}{l}v^{(\alpha)}=\xi_{\pos_l}^{(\alpha)},\\
\alpha=1,\dots,m\end{array}$}}}\!\!;\]
As $E_{{\cal{W}}}E_{{\cal{U}}}=E_{{\cal{U}}}$, $E_{{\cal{U}},\pos_n}^2=E_{{\cal{U}},\pos_n}$ and $\pos_l=\pos_n$, we have
\begin{equation}\label{calc1}\begin{array}{ll}E_{{\cal{W}},\pos_l}s_ls_{l+1}\dots s_{n-1}E_{{\cal{U}},\pos_n}& = \displaystyle s_ls_{l+1}\dots s_{n-1}\cdot \prod_{\alpha=1}^m\ \frac{v^{(\alpha)}-\xi_{\pos_l}^{(\alpha)}}{v^{(\alpha)}-g^{(\alpha)}_{n}}E_{{\cal{W}}} E_{{\cal{U}},\pos_n}\Bigr\vert_{\hspace{-0.15cm}\textrm{\scriptsize{$\begin{array}{l}v^{(\alpha)}=\xi_{\pos_l}^{(\alpha)},\\
\alpha=1,\dots,m\end{array}$}}}\\[2em]
 & = s_ls_{l+1}\dots s_{n-1}E_{{\cal{U}},\pos_n}\ .\end{array}\end{equation}
Thus, we rewrite (\ref{lhs1}) as
\[F_{_{{\cal{T}}}}(\di_{\pos_n}u)\,s_{n-1}\dots s_{l+1}\bigl(s_l+\frac{e_{l,l+1}}{u-\cc^G_l}\bigr)\cdot\phi_l(\cc^G_1,\dots,\cc^G_{l-1},u) E_{{\cal{W}},\pos_l}\cdot s_ls_{l+1}\dots s_{n-1} E_{{\cal{U}},\pos_n}\ .\]
We use the induction hypothesis to replace $\phi_l(\cc^G_1,\dots,\cc^G_{l-1},u) E_{{\cal{W}},\pos_l}$ by $\displaystyle \bigl(F_{_{{\cal{V}}}}(\di_{\pos_n}u)\bigr)^{-1}\frac{u-\cc^G_{l}}{u-\tj_{l}}\ E_{{\cal{W}},\pos_l}$, and we use (\ref{calc1}) again to obtain for the left hand side of (\ref{eq-idem12}):
\[F_{_{{\cal{T}}}}(\di_{\pos_n}u)(F_{_{{\cal{V}}}}(\di_{\pos_n}u))^{-1}\,s_{n-1}\dots s_{l+1}\bigl(s_l+\frac{e_{l,l+1}}{u-\cc^G_l}\bigr)\,\frac{u-\cc^G_{l}}{u-\tj_{l}}\,s_ls_{l+1}\dots s_{n-1} E_{{\cal{U}},\pos_n}\ .\]
Recall that the inverse of $\displaystyle\ \bigl(s_l+\frac{e_{l,l+1}}{u-\cc^G_l}\bigr)\ $ is $\ \displaystyle \bigl(s_l+\frac{e_{l,l+1}}{\cc^G_l-u}\bigr)\Bigl(1-\frac{e_{l,l+1}^2}{(u-\cc^G_{l})^2}\Bigr)^{-1}\ $ according to the second line in (\ref{YB2}). We move $\displaystyle\ s_{n-1}\dots s_{l+1}\bigl(s_l+\frac{e_{l,l+1}}{u-\cc^G_l}\bigr)\,(u-\tj_{l})^{-1}\,$ to the right hand side of (\ref{eq-idem12}) and, using that $\tj_n$ commutes with $E_{{\cal{U}},\pos_n}$, we move $(u-\tj^{\phantom{\text{\Large{A}}}}_n\!\!\!)^{-1}$ from the right hand side of (\ref{eq-idem12}) to the left hand side. We finally obtain that$^{\phantom{\text{\large{A}}}}\!\!\!\!\!$ the equality (\ref{eq-idem12}) is equivalent to:
\begin{equation}\label{eq-inter11}
\begin{array}{l}\displaystyle F_{_{{\cal{T}}}}(\di_{\pos_n}u)(F_{_{{\cal{V}}}}(\di_{\pos_n}u))^{-1}(u-\cc^G_{l}) s_ls_{l+1}\dots s_{n-1}(u-\tj_n) \, E_{{\cal{U}},\pos_n}\\[1em]
\hspace{2cm}=\displaystyle(u-\cc^G_n)(u-\tj_{l})\bigl(s_l+\frac{e_{l,l+1}}{\cc^G_l-u}\bigr)\Bigl(1-\frac{e_{l,l+1}^2}{(u-\cc^G_{l})^2}\Bigr)^{-1}\! s_{l+1}\dots s_{n-1}\, E_{{\cal{U}},\pos_n}\ ;\end{array}\end{equation}
Now notice that $e_{l,l+1}\cdot s_{l+1}\dots s_{n-1}=s_{l+1}\dots s_{n-1}\cdot e_{l,n}$ and that (see (\ref{eq-eij-rep})) $e_{l,n}^2\,E_{{\cal{U}},\pos_n}=\displaystyle\frac{1}{\di_{\pos_n}^2}\,E_{{\cal{U}},\pos_n}$ since $\pos_l=\pos_n$.
Moreover, formula (\ref{def-F1}) implies, since $\pos_l=\pos_n$ and $\pos_i\neq\pos_n$ for $i=l+1,\dots,n-1$, 
\begin{equation}\label{coeff-inter11}F_{_{{\cal{T}}}}(\di_{\pos_n}u)(F_{_{{\cal{V}}}}(\di_{\pos_n}u))^{-1} = \,\frac{\di_{\pos_n}u-\cc_{n}}{\di_{\pos_n}u-\cc_l}\frac{(\di_{\pos_n}u-\cc_l)^2}{(\di_{\pos_n}u-\cc_l)^2-1}\,=\,\frac{u-\cc^G_{n}}{u-\cc^G_l}\frac{(u-\cc^G_l)^2}{(u-\cc^G_l)^2-\frac{1}{\di^2_{\pos_n}}}.\end{equation}
Therefore, to verify formula (\ref{eq-inter11}), it remains to show that
\begin{equation}\label{eq-inter21}
s_ls_{l+1}\dots s_{n-1}(u-\tj_n) \, E_{{\cal{U}},\pos_n}\,=\,(u-\tj_{l})\bigl(s_l+\frac{e_{l,l+1}}{\cc^G_l-u}\bigr) s_{l+1}\dots s_{n-1}\, E_{{\cal{U}},\pos_n}\ .
\end{equation}
Rewrite the right hand side of (\ref{eq-inter21}) as $$(u-\tj_{l})s_ls_{l+1}\dots s_{n-1}\, E_{{\cal{U}},\pos_n}+ (u-\tj_l)\frac{e_{l,l+1}}{\cc^G_l-u} s_{l+1}\dots s_{n-1}\, E_{{\cal{U}},\pos_n}\ .$$
Then notice that $$\tj_l\,e_{l,l+1}\,s_{l+1}\dots s_{n-1}\, E_{{\cal{U}},\pos_n}=s_{l+1}\dots s_{n-1}\tj_l\,e_{l,n}\, E_{{\cal{U}},\pos_n}\ .$$
We use that $\tj_l\,e_{l,n}\, E_{{\cal{U}},\pos_n}=\tj_l\, E_{{\cal{U}},\pos_n}\,e_{l,n}=\cc^G_l\,e_{l,n}\, E_{{\cal{U}},\pos_n}$ and the Lemma \ref{lemm-tj} to replace $\tj_{l}s_ls_{l+1}\dots s_{n-1}$; we obtain for the right hand side of (\ref{eq-inter21})
\[ \Bigl( s_ls_{l+1}\dots s_{n-1}(u-\tj_n)+\sum_{k=l}^{n-1}s_ls_{l+1}\dots \widehat{s_k}\dots s_{n-1}\,e_{k,n}-s_{l+1}\dots s_{n-1}e_{l,n} \Bigr)\, E_{{\cal{U}},\pos_n}\ .\]
As $e_{k,n}\,E_{{\cal{U}},\pos_n}=0$ if $k>l$, the formula (\ref{eq-inter21}) is verified.\hfill $\square$

\vskip .2cm
\textbf{Remark.} During the proof of the Theorem \ref{prop-fus1}, we transformed the defining formula (\ref{def-Phi1}) for the rational function $\Phi(u_1,\dots,u_n,\underline{v_1},\dots,\underline{v_n})$ into the formula (\ref{Phi-ind}). The formula (\ref{Phi-ind}) is well adapted to the structure of chain with respect to $n$ of the groups $\tG_n$, in the sense that, using formula (\ref{Phi-ind}), we write 
$$\Phi(u_1,\dots,u_n,\underline{v_1},\dots,\underline{v_n})=\widetilde{\phi}_n(u_1,\dots,u_n,\underline{v_n})\Phi^{(n-1)}(u_1,\dots,u_{n-1},\underline{v_1},\dots,\underline{v_{n-1}})\ ,$$
where $\Phi^{(n-1)}(u_1,\dots,u_{n-1},\underline{v_1},\dots,\underline{v_{n-1}})$ is the rational function corresponding to $\tG_{n-1}$ (seen as a rational function with values in $\tG_n$) and $\widetilde{\phi}_n(u_1,\dots,u_n,\underline{v_n})$ is defined by (\ref{def-phit}). The formula (\ref{Phi-ind}) is often useful for explicit calculations.

\paragraph{Examples.} Here we consider the example when $G$ is the symmetric group on $3$ letters. To avoid confusion we will keep the notation $G$ for this group, and $\tG_n$ for the wreath product of $G$ by the symmetric group $S_n$. We use the standard cyclic notation for permutations and we denote the elements of $G$ by $1_{G}$, $(1,2)$, $(1,3)$, $(2,3)$, $(1,2,3)$ and $(1,3,2)$. The central elements, defined in (\ref{def-g}), of $\C G$ are 
$$g^{(1)}:=\frac{1}{3}\Bigl( (1,2)+(1,3)+(2,3)\Bigr),\quad g^{(2)}:=\frac{1}{2}\Bigl( (1,2,3)+(1,3,2)\Bigr)\quad\text{and}\quad g^{(3)}:=1_{G}\ .$$
Let $\rho_1,\rho_2,\rho_3$ be the pairwise non-isomorphic irreducible representations of $G$, namely $\rho_1$ is the trivial representation, $\rho_2$ is the sign representation and $\rho_3$ is the two-dimensional irreducible representation of $G$. The numbers defined in (\ref{def-p}) are equal to:
$$\xi_1^{(1)}=1,\ \xi_2^{(1)}=-1,\ \xi_3^{(1)}=0,\ \ \ \ \xi_1^{(2)}=1,\ \xi_2^{(2)}=1,\ \xi_3^{(2)}=-\frac{1}{2}\ \ \ \ \text{and}\ \ \ \ \xi_1^{(3)}=1,\ \xi_2^{(3)}=1,\ \xi_3^{(3)}=1\ .$$
Thus the functions defined by (\ref{bax-g}) are, in this example,
$$g^{(1)}(v)=\bigl(g^{(1)}\bigr)^2+v g^{(1)}+v^2-1\,,\ \ \ \  g^{(2)}(v)=g^{(2)}+v-\frac{1}{2}\ \ \ \ \text{and}\ \ \ \ g^{(3)}(v)=1\ .$$
Recall the notation $g^{(\alpha)}_i$, $\alpha=1,2,3$ and $i=1,\dots,n$, for images of $g^{(\alpha)}$ in $\C\tG_n$ (see (\ref{inj-Gn})).

\vskip .2cm
For the 3-partition of size 1 $\left(\Box,\varnothing,\varnothing\right)$, we have, from formulas (\ref{def-fg}) and (\ref{def-f}), $\mathsf{F}^G_{\left(\Box,\varnothing,\varnothing\right)}=3$ and $\mathsf{F}_{\left(\Box,\varnothing,\varnothing\right)}=1$.

We obtain for the idempotent $E_{\left(\fbox{\scriptsize{$1$}},\,\varnothing,\,\varnothing\right)}$ of $\C\tG_1$($\cong \C G$):
$$E_{\left(\fbox{\scriptsize{$1$}},\,\varnothing,\,\varnothing\right)}=\frac{1}{3}\Bigl(\bigl(g_1^{(1)}\bigr)^2+g_1^{(1)}\Bigr)\Bigl(g_1^{(2)}+\frac{1}{2}\Bigr)\ .$$

For the 3-partition of size 2 $\left(\Box,\varnothing,\Box\right)$, we have $\mathsf{F}^G_{\left(\Box,\varnothing,\Box\right)}=\frac{9}{2}$ and $\mathsf{F}_{\left(\Box,\varnothing,\Box\right)}=1$. We obtain, from Theorem \ref{prop-fus1} and using formula (\ref{Phi-ind}),
$$E_{\left(\fbox{\scriptsize{$1$}},\,\varnothing,\,\fbox{\scriptsize{$2$}}\right)}=\frac{2}{9}\,s_1(\cc_2,0)\,s_1\cdot\Bigl(\bigl(g_2^{(1)}\bigr)^2-1\Bigr)\Bigl(g_2^{(2)}-1\Bigr)\cdot 3E_{\left(\fbox{\scriptsize{$1$}},\,\varnothing,\,\varnothing\right)}\Biggr\vert_{\cc_2=0}\ .$$
In particular, Theorem \ref{prop-fus1} asserts that the above function is non-singular at $\cc_2=0$ although, from (\ref{def-baxt}), $s_1(\cc_2,0)$ is singular at $\cc_2=0$. As shows the proof of the Lemma \ref{lemm-idem11}, this comes from the fact that 
$$e_{1,2}\cdot\Bigl(\bigl(g_2^{(1)}\bigr)^2-1\Bigr)\Bigl(g_2^{(2)}-1\Bigr)\cdot E_{\left(\fbox{\scriptsize{$1$}},\,\varnothing,\,\varnothing\right)}=0\ .$$
So actually we have
$$\begin{array}{cl} E_{\left(\fbox{\scriptsize{$1$}},\,\varnothing,\,\fbox{\scriptsize{$2$}}\right)} & =\ \displaystyle\frac{2}{9}\,\Bigl(\bigl(g_2^{(1)}\bigr)^2-1\Bigr)\Bigl(g_2^{(2)}-1\Bigr)\cdot 3E_{\left(\fbox{\scriptsize{$1$}},\,\varnothing,\,\varnothing\right)}\\[1em]
 & =\ \displaystyle\frac{2}{9}\,\Bigl(\bigl(g_2^{(1)}\bigr)^2-1\Bigr)\Bigl(g_2^{(2)}-1\Bigr)\Bigl(\bigl(g_1^{(1)}\bigr)^2+g_1^{(1)}\Bigr)\Bigl(g_1^{(2)}+\frac{1}{2}\Bigr)\ .\end{array}$$

For the 3-partition of size 3 $\left(\Box\!\Box,\varnothing,\Box\right)$, we have $\mathsf{F}^G_{\left(\Box\!\Box,\varnothing,\Box\right)}=\frac{27}{2}$ and $\mathsf{F}_{\left(\Box\!\Box,\varnothing,\Box\right)}=2$. We obtain, from Theorem \ref{prop-fus1} and using formula (\ref{Phi-ind}),
$$E_{\left(\fbox{\scriptsize{$1$}}\fbox{\scriptsize{$3$}},\,\varnothing,\,\fbox{\scriptsize{$2$}}\right)}=\ \frac{1}{27}\,s_2(1,0)s_1(1,0)\,s_1s_2\,\Bigl(\bigl(g_3^{(1)}\bigr)^2+g_3^{(1)}\Bigr)\Bigl(g_3^{(2)}+\frac{1}{2}\Bigr)\cdot \frac{9}{2}E_{\left(\fbox{\scriptsize{$1$}},\,\varnothing,\,\fbox{\scriptsize{$2$}}\right)}\ .$$
Note that, as shows the proof of the Lemma \ref{lemm-idem11}, $e_{2,3}\cdot s_1(1,0)\,s_1s_2\,\Bigl(\bigl(g_3^{(1)}\bigr)^2+g_3^{(1)}\Bigr)\Bigl(g_3^{(2)}+\frac{1}{2}\Bigr)\cdot E_{\left(\fbox{\scriptsize{$1$}},\,\varnothing,\,\fbox{\scriptsize{$2$}}\right)}=0$, and so we have actually
 $$E_{\left(\fbox{\scriptsize{$1$}}\fbox{\scriptsize{$3$}},\,\varnothing,\,\fbox{\scriptsize{$2$}}\right)}=\ \frac{1}{27}\,s_2s_1(1,0)\,s_1s_2\,\Bigl(\bigl(g_3^{(1)}\bigr)^2+g_3^{(1)}\Bigr)\Bigl(g_3^{(2)}+\frac{1}{2}\Bigr)\cdot \frac{9}{2}E_{\left(\fbox{\scriptsize{$1$}},\,\varnothing,\,\fbox{\scriptsize{$2$}}\right)}\ .$$

\vskip .2cm
\textbf{Remark.} We notice that the central element $g^{(3)}=1_{G}$ does not appear in the fusion formula (because $\xi_1^{(3)}=\xi_2^{(3)}=\xi_3^{(3)}$). Moreover, in the examples above, we could have used only the central element $g^{(1)}$ since its eigenvalues are enough to distinguish between the irreducible representations of $G$; namely, we would have obtained
$$E_{\left(\fbox{\scriptsize{$1$}},\,\varnothing,\,\varnothing\right)}=\frac{1}{2}\Bigl(\bigl(g_1^{(1)}\bigr)^2+g_1^{(1)}\Bigr)\ ,$$
$$E_{\left(\fbox{\scriptsize{$1$}},\,\varnothing,\,\fbox{\scriptsize{$2$}}\right)} =\ -\frac{1}{2}\,\Bigl(\bigl(g_2^{(1)}\bigr)^2-1\Bigr)\Bigl(\bigl(g_1^{(1)}\bigr)^2+g_1^{(1)}\Bigr)\ ,$$  
$$E_{\left(\fbox{\scriptsize{$1$}}\fbox{\scriptsize{$3$}},\,\varnothing,\,\fbox{\scriptsize{$2$}}\right)}=\ -\frac{1}{8}\,s_2s_1(1,0)\,s_1s_2\,\Bigl(\bigl(g_3^{(1)}\bigr)^2+g_3^{(1)}\Bigr) \Bigl(\bigl(g_2^{(1)}\bigr)^2-1\Bigr)\Bigl(\bigl(g_1^{(1)}\bigr)^2+g_1^{(1)}\Bigr)\ .$$

In general, the central element corresponding to the conjugacy class of the unit element of $G$ never appears. Nevertheless, in general, all other central elements (\ref{def-g}) of $\C G$ are necessary. In the following Section, we explain how to reduce the number of elements $g^{(\alpha)}$ appearing in the fusion procedure for the particular situation of an Abelian finite group $G$. 

\section{{\hspace{-0.55cm}.\hspace{0.55cm}}Simplified fusion formula when $G$ is Abelian}\label{sec-abe}

From now let $G$ be any finite Abelian group. It is a standard fact that $G$ is thus isomorphic to a direct product of finite cyclic groups. So we assume that
\begin{equation}\label{G-abel}G\cong C_{k_1}\times C_{k_2}\times\dots\times C_{k_N}\ ,\end{equation}
where $N$ is a non-negative integer, $k_1,k_2,\dots,k_N\geq2$ and $C_{k_{\alpha}}$ is the cyclic group of order $k_{\alpha}$, $\alpha=1,\dots,N$. Set $m:=k_1k_2\dots k_N$, the cardinality of $G$.

\vskip .2cm
As $G$ is Abelian, the number of its conjugacy classes is equal to $m$. The set of central elements $g^{(\alpha)}$ of $\C G$ defined in (\ref{def-g}) coincide here with the set of elements of $G$. All these elements appear in the fusion formula of the previous Section, see (\ref{def-gamma}). We will provide a simplified formula using a certain subset of elements of $G$.

\vskip .2cm
For $\alpha=1,\dots,N$, we choose and denote by $t^{(\alpha)}$ an element of $G$ which generates the cyclic group $C_{k_{\alpha}}$ appearing in (\ref{G-abel}). Notice that the elements $t^{(\alpha)}$, $\alpha=1,\dots,N$, satisfy the relations
\[\bigl(t^{(\alpha)}\bigr)^{k_{\alpha}}=1\ \ \ \ \text{and}\ \ \ \ t^{(\alpha)}t^{(\alpha')}=t^{(\alpha')}t^{(\alpha)}\qquad\ \ \text{for $\alpha,\alpha'=1,\dots,N$.}\]
Moreover the elements $t^{(\alpha)}$, $\alpha=1,\dots,N$, generate the group $G$ and thus the knowledge of their eigenvalues is sufficient to distinguish between the irreducible representations of $G$. Namely, any irreducible representation of $G$ is one-dimensional and is obtained by sending $t^{(\alpha)}$ to a $k_{\alpha}$-th root of unity for $\alpha=1,\dots,N$ (there are $m$ non-isomorphic representations of this sort and they exhaust the set of irreducible representations of $G$).

\vskip .2cm
So now let $S^{(\alpha)}:=\{\xi_1^{(\alpha)},\dots,\xi_{k_{\alpha}}^{(\alpha)}\}$ be the set of all $k_{\alpha}$-th roots of unity, and, for $i=1,\dots,n$, define as in (\ref{bax-gi})
\begin{equation}\label{bax-gi'}
t^{(\alpha)}_i(v):=\frac{\prod_{\xi^{(\alpha)}\in S^{(\alpha)}}(v-\xi^{(\alpha)})}{v-t_i^{(\alpha)}}\ \ \ \ \ \text{for $\alpha=1,\dots,N$\,,}
\end{equation}
where we recall that $t^{(\alpha)}_i$, $i=1,\dots,n$, is the image of the element $t^{(\alpha)}$ by the injective morphism $\iota\circ\varpi_i$ from $\C G$ to $\C\tG_n$, see (\ref{inj-Gn}).

\vskip .2cm
From (\ref{bax-g2}), it is straightforward to see that we have, for $i=1,\dots,n$ and $\alpha=1,\dots,N$,
\[t_i^{(\alpha)}(v)=v^{k_{\alpha}-1}+v^{k_{\alpha}-2}t_i^{(\alpha)}+\dots+v\,\bigl(t_i^{(\alpha)}\bigr)^{k_{\alpha}-2}+\bigl(t_i^{(\alpha)}\bigr)^{k_{\alpha}-1}\ .\]

We define
\begin{equation}\label{def-gamma'}\Gamma'(\underline{v_1},\dots,\underline{v_n}):=\displaystyle\prod_{i=1}^n\prod_{\alpha=1}^N t^{(\alpha)}_i(v^{(\alpha)}_i)\ ,\end{equation}
and, for an $m$-partition $\lam$  (compare with (\ref{def-fg})),
\begin{equation}\label{def-fg2}
\mathsf{F'}^G_{\lam}:=\prod_{\thm\in\lam}
\Biggl(\ \prod_{\alpha=1}^{N}\Bigl(\!\!\prod_{\textrm{\scriptsize{$\begin{array}{c}\xi^{(\alpha)}\in S^{(\alpha)}\\[0em]
\xi^{(\alpha)}\neq\xi_{\pos(\thm)}^{(\alpha)}\end{array}$}}
}\!\!\!\!\bigl(\xi_{\pos(\thm)}^{(\alpha)}-\xi^{(\alpha)}\bigr)\ \Bigr)\,\Biggr)
\ ,
\end{equation}
The formulation of the Theorem \ref{prop-fus1} remains the same with $\Gamma(\underline{v_1},\dots,\underline{v_n})$ replaced by $\Gamma'(\underline{v_1},\dots,\underline{v_n})$ in the defining formula (\ref{def-Phi1}) for the function $\Phi(u_1,\dots,u_n,\underline{v_1},\dots,\underline{v_n})$\,, and $\mathsf{F}^G_{\lam}$ replaced by $\mathsf{F'}^G_{\lam}$.

\vskip .2cm
\textbf{Remark.} Let $\zeta$ be a $k$-th root of unity. Recall that $\prod_{\xi} (\zeta-\xi)=k\zeta^{-1}$, where the product is taken on the $k$-th roots of unity $\xi$ different from $\zeta$. Thus, for $G$ as in (\ref{G-abel}), we have for the coefficients $\mathsf{F'}^G_{\lam}$ defined by (\ref{def-fg2}):
\[\mathsf{F'}^G_{\lam}=\prod_{\thm\in\lam}\Bigl(\ \prod_{\alpha=1}^{N}\frac{k_{\alpha}}{\xi^{(\alpha)}_{\pos(\thm)}}\Bigr)\ .\]

\paragraph{Example.} Consider the situation $N=1$ and $k_1=3$ in (\ref{G-abel}). Thus $G$ is the cyclic group of order 3, and $\tG_n$ is the complex reflection group of type $G(3,1,n)$. Let $t:=t^{(1)}$ and let $\{\xi_1,\xi_2,\xi_3\}$ be the set of all third roots of unity. The Theorem \ref{prop-fus1} (in its simplified version explained above) asserts that, for example, the idempotent $E_{\left(\fbox{\scriptsize{$1$}}\fbox{\scriptsize{$3$}},\,\varnothing,\,\fbox{\scriptsize{$2$}}\right)}$ can be expressed as
$$E_{\left(\fbox{\scriptsize{$1$}}\fbox{\scriptsize{$3$}},\,\varnothing,\,\fbox{\scriptsize{$2$}}\right)}=\ \frac{\xi_1^2\xi_3}{54}s_2(1,0)s_1(1,0)s_1s_2s_1(0,0)s_1\,(\xi_1^2+\xi_1t_3+t_3^2)(\xi_3^2+\xi_3t_2+t_2^2)(\xi_1^2+\xi_1t_1+t_1^2)\,.$$
As already seen in the examples of the preceding Section, the proof of the Theorem \ref{prop-fus1} shows that we actually have:
$$E_{\left(\fbox{\scriptsize{$1$}}\fbox{\scriptsize{$3$}},\,\varnothing,\,\fbox{\scriptsize{$2$}}\right)}=\ \frac{\xi_1^2\xi_3}{54}s_2s_1(1,0)s_1s_2\,(\xi_1^2+\xi_1t_3+t_3^2)(\xi_3^2+\xi_3t_2+t_2^2)(\xi_1^2+\xi_1t_1+t_1^2)\,.$$

\end{document}